\documentclass{amsart}

\newtheorem{thm}{Theorem}
\newtheorem*{thma}{Theorem 5(b)}

\newtheorem{lem}[thm]{Lemma}
\newtheorem{prop}[thm]{Proposition}
\newtheorem{prob}[thm]{Problem}
\theoremstyle{definition}
\newtheorem{defn}[thm]{Definition}
\theoremstyle{remark}

\numberwithin{equation}{section}

\newcommand{\To}{\longrightarrow}

\begin{document}

\title{Some examples of continuous images of Radon-Nikod\'ym compact spaces}

\author{Alexander D. Arvanitakis, Antonio Avil\'{e}s}

\thanks{The second author was supported by a Marie Curie
Intra-European Fellowship MCEIF-CT2006-038768 and research
projects MTM2005-08379 (MEC and FEDER) and S\'{e}neca 00690/PI/04}

\address{National Technical University of
Athens, Department of Mathematics, Athens 15780, Greece}
\email{aarva@math.ntua.gr}

\address{Universit\'e de Paris VII, Equipe de Logique
Math\'ematique, UFR de Math\'ematiques, 2 Place Jussieu, 75251
Paris, France} \email{avileslo@um.es, aviles@logique.jussieu.fr}.

\begin{abstract} We provide a characterization of continuous images
of Radon-Nikod\'ym compacta lying in a product of real lines and
model on it a method for constructing natural examples of such
continuous images.\end{abstract}



\maketitle

\section{Introduction}\label{sec1}

This note contains some ideas to address the well known open
problem whether the class of Radon-Nikod\'ym (or briefly RN)
compact spaces is stable under continuous images. In particular,
we provide a method to construct a number of concrete compact
spaces which seem to be natural candidates to be counterexamples,
but we are unable to decide whether they are Radon-Nikod\'ym
compact or not. We think that the understanding of these examples
could be a useful tool towards the solution of the problem.

\section{Basic definitions and notations}

We denote by $\mathcal{N}=\mathbb{N}^\mathbb{N}$ the set of all sequences of natural numbers, and we consider that $0\in\mathbb{N}$. If $s\in\mathbb{N}^{<\omega}$ is a finite sequence of natural numbers and $\sigma\in\mathcal{N}$, then $s<\sigma$ means that $s$ is an initial segment of $\sigma$. We call $\mathcal{N}_s = \{ \sigma\in \mathcal{N} : s<\sigma\}$.\\

If $K$ is a topological space, a map $d:K\times K\To \mathbb{R}$
is said to $\varepsilon$-fragment $K$ if for every (closed) set
$L\subset K$ there exists a nonempty relative open subset $U$ of
$L$ of $d$-diameter less than $\varepsilon$, that is,
$\sup\{d(x,y) : x,y\in U\}<\varepsilon$. If $d$
$\varepsilon$-fragments $K$ for every $\varepsilon$ we shall say
that $d$ fragments $K$.\\

Given a bounded family $\Delta$ of continuous functions over $K$
we shall denote by $d_\Delta$ the uniform pseudometric over
$\Delta$, that is $d_\Delta(x,y) =\sup\{|f(x)-f(y)| :
f\in\Delta\}$. When we view a compactum lying as a subset of a
product of real lines, $K\subset \mathbb{R}^\Gamma$, we shall
identify every element $\gamma\in\Gamma$ with the continuous
function on $K$ given by projection onto the $\gamma$-th
coordinate. Thus, for a given $\Delta\subset\Gamma$,
$d_\Delta(x,y)
= \sup\{|x_\gamma-y_\gamma| : \gamma\in\Delta\}$.\\

Originally, RN compact spaces are defined as weak$^\ast$ compact subsets of dual Banach spaces with the Radon-Nikod\'ym property. From an intrinsic topological point of view, they can be characterized as those compact spaces $K$ for which there exists a lower semicontinuous metric $d:K\times K\To \mathbb{R}$ which fragments $K$. The class of quasi Radon-Nikod\'ym (qRN) compact spaces is a superclass of the class of RN compacta which is stable under continuous images. According to the definition of the first author \cite{Arvanitakis}, a compact space is qRN compact if there exists a lower semicontinuous map $f:K\times K\To\mathbb{R}^+$ which is nonzero out of the diagonal and which fragments $K$. This is actually equivalent to other definitions given by Fabian, Heiler and Matouskova~\cite{FHM} and Reznichenko, as it is shown in \cite{Namiokanote} and \cite{cardinalb}. We refer to the recent survey of Fabian~\cite{Fabiansurvey} on this subject, as well as other articles containing different approaches or partial results concerning the problem of continuous images as \cite{OSV}, \cite{IW} or \cite{Aviorder}.\\

In this note, we will work with compact spaces viewed as closed subsets of Tychonoff cubes $[-1,1]^\Gamma$. We present below the two known characterizations of both RN and qRN compact spaces in terms of these embeddings. Theorem~\ref{RN} can be found in \cite{Namioka}. With respect to Theorem~\ref{qRN}, the equivalence of $(2)$ and $(3)$ is proven in \cite{FHM}, and the equivalence of $(1)$ with the others in \cite{cardinalb}.\\

\begin{thm}\label{RN}
A compact space $K$ is RN compact if and only if there exists
an embedding of $K$ into a product of intervals $K\subset
[-1,1]^\Gamma$ such
that $d_\Gamma$ fragments $K$.\\
\end{thm}

\begin{thm}\label{qRN}
For a compact space $K$ the following are equivalent:
\begin{enumerate}\item $K$ is quasi RN compact. \item There exists
an embedding of $K$ into a product of intervals $K\subset
[-1,1]^\Gamma$ and a map $u:\Gamma\To\mathcal{N}$ such that for
every finite sequence $s$ of naturals of length $n>0$ the
pseudometric $d_{u^{-1}(\mathcal{N}_s)}$ $\frac{1}{2^n}$-fragments
$K$. \item For any embedding of $K$ into a product of intervals
$K\subset [-1,1]^\Gamma$ there exists a map
$u:\Gamma\To\mathcal{N}$ such that for every finite sequence $s$
of naturals of length $n>0$ the pseudometric
$d_{u^{-1}(\mathcal{N}_s)}$  $\frac{1}{2^n}$-fragments
$K$.\\
\end{enumerate}
\end{thm}

The picture is that we have three classes of compact sets, RN
compacta, their continuous images (ciRN) and quasi RN (qRN)
compacta,

$$ RN \subset ciRN \subset qRN$$

and we do not know at all whether any of the inclusions is
strict.\\

A key conceptual difference between RN and qRN is that there is no
analogue of condition (3) in Theorem~\ref{qRN} for RN, referring
to \emph{any embedding} of $K$ into a cube. This makes the class
qRN easier to handle than RN in many aspects: in order to check
that some space $K$ is qRN compact we take any embedding of $K$
and we try to see whether there is function
$u:\Gamma\To\mathcal{N}$ fulfilling the mentioned condition (3).
However, if we want to check that a certain space is RN compact we
must find an appropriate embedding, probably different from the
obvious ones. Such a difficulty is
found, of course, in the problem of the continuous image.\\

\section{A characterization of continuous images of RN compacta}

The main result of this note, and the inspiration for the
announced concrete examples, will be Theorem~\ref{ciRN}, where we
give a similar characterization as those appearing in
Theorems~\ref{RN} and~\ref{qRN} for the intermediate class ciRN of
the continuous images of Radon-Nikod\'ym compacta.\\

\begin{defn}
A family $\Delta\subset C(K,[-1,1])$ of continuous functions over
$K$ is said to be a Namioka family if the pseudometric $d_\Delta$
fragments $K$.\\
\end{defn}

Using this concept, Theorem~\ref{RN} can be restated saying that a
compact space is RN if and only if there exists a Namioka family
$\Delta\subset C(K)$ which separates the points of $K$. We shall
need the fact that if a Namioka family exists, then indeed it can
be chosen to be much bigger than simply a separating
family.\\

\begin{lem}\label{lemNamiokafamily}
Let $K$ be an RN compactum. Then there exists a Namioka family
$\Delta$ such that $C(K)
=\overline{\bigcup_{n\in\mathbb{N}}n\Delta}^{\|\cdot\|}$.\\
\end{lem}

PROOF: We consider $\Delta_0$ some Namioka family over $K$ and
then we define $\Delta$ to be the set of all $f\in C(K)$ such that
$|f(x)-f(y)|\leq d_{\Delta_0}(x,y)$ for all $x,y\in K$. Clearly
$d_\Delta = d_{\Delta_0}$ so $\Delta$ is a Namioka family. On the
other hand, $\bigcup_{n\in\mathbb{N}}n\Delta $ is a linear lattice
of functions which separates the points of $K$, so by the
Stone-Weierstrass Theorem, it is uniformly dense in $C(K)$.$\qed$\\

\begin{thm}\label{ciRN}
For a compact space $K\subset[a,b]^\Gamma$ the following are
equivalent:\begin{enumerate} \item $K$ is the continuous image of
an RN compactum. \item There exists a function
$u:\Gamma\To\mathcal{N}$ and a family of compact sets $K_s\subset
[a,b]^\Gamma$ for $s\in\mathbb{N}^{<\omega}$ such
that\begin{enumerate} \item $d_{u^{-1}(\mathcal{N}_s)}$ fragments
$K_s$ for every $s$, and \item If $length(s)=n>0$, then for every
$x\in K$ there exists $y\in K_s$ such that
$d_{u^{-1}(\mathcal{N}_s)}(x,y)\leq\frac{1}{2^n}$.\\
\end{enumerate}
\end{enumerate}
\end{thm}

Of course, the use of the numbers $\frac{1}{2^n}$ is inessential,
we could have used any numbers $\varepsilon_n$ converging to 0. We
give an equivalent reformulation of Theorem~\ref{ciRN} which will
be more suitable for the purpose of making the proof more
transparent, though it will be the previous statement which we
will be more relevant in further discussion.\\

\begin{thma}
For a compact space $K\subset [a,b]^\Gamma$ the following are
equivalent:\begin{enumerate} \item $K$ is a continuous image of an
RN compactum. \item For every $\varepsilon>0$ there exists a
countable decomposition $\Gamma=\bigcup_m\Gamma^\varepsilon_m$ and
compact sets $K^\varepsilon_m\subset[a,b]^\Gamma$ such
that:\begin{enumerate} \item $K_m^\varepsilon$ is fragmented by
$d_{\Gamma^\varepsilon_m}$ \item For every $m\in\mathbb{N}$ and
$x\in K$ there is $y\in K^\varepsilon_m$
such that $d_{\Gamma^\varepsilon_m}(x,y)\leq\varepsilon$.\\
\end{enumerate}
\end{enumerate}
\end{thma}

Before passing to the proof, we state a lemma from \cite{Arvanitakis} that we will often use.\\

\begin{lem}\label{quotientfrag}
Let $f:Q\To S$ be a continuous surjection between compact spaces, and $d:Q\times Q\To \mathbb{R}$ a lower semicontinuous map that fragments $Q$. Then the map $\hat{d}(x,y) = \inf\{d(u,v) : f(u)=x, f(v)=y\}$ fragments $S$.\\
\end{lem}

\emph{Proof}: Theorem~\ref{ciRN} is easily seen to be equivalent
to Theorem~\ref{ciRN}(\textbf{b}), we leave this to the reader. We
shall prove Theorem~\ref{ciRN}(\textbf{b}). Without loss of
generality,
we suppose that $[a,b] = [-2,2]$.\\

$(1) \Rightarrow (2)$: Let $L$ be RN compact and $\pi:L\To
K\subset[-2,2]^\Gamma$ a continuous surjection. We take a Namioka
family $\Delta$ on $L$ as in Lemma~\ref{lemNamiokafamily}, so that
we view $L\subset [-1,1]^\Delta$. We fix $\varepsilon\in (0,1)$.
For every $\gamma\in\Gamma$ we have a continuous function on $L$
given by $z\mapsto \pi(z)_\gamma$, so since $C(L)
=\overline{\bigcup_{m\in\mathbb{N}}m\Delta}^{\|\cdot\|}$ we have
that:\\

\emph{for every $\gamma\in\Gamma$ there exists
$m(\gamma)\in\mathbb{N}$ and $\delta(\gamma)\in\Delta$ such for
every $z\in L$,
$|\pi(z)_\gamma - m(\gamma)z_{\delta(\gamma)}|\leq\varepsilon$.}\\

 We set
$$\Gamma_m^\varepsilon = \{\gamma : m(\gamma)=m\}$$
$$K_m^\varepsilon = \{x\in [-2,2]^\Gamma : \exists z\in L :
m z_{\delta(\gamma)} = x_\gamma
\forall\gamma\in\Gamma_m^\varepsilon\}$$

Each $K_m^\varepsilon$ is compact, because it is the continuous image of $L$ under the map $z\mapsto (mz_{\delta(\gamma)})_{\gamma\in\Gamma}$.\\

 We check condition (b). If $x\in K$, then $x=\pi(z)$ for some $z\in L$, and we know that
$|x_\gamma - m(\gamma)z_{\delta(\gamma)}|\leq\varepsilon$ for
every $\gamma\in\Gamma$. We define $y$ to be an element of
$[-2,2]^\Gamma$ such that $y_\gamma = m z_{\delta(\gamma)}$ for
$\gamma\in\Gamma_m^\varepsilon$. Then $y\in K_m^\varepsilon$ and
$d_{\Gamma_m^\varepsilon}(x,y)\leq\varepsilon$.\\

 For condition (a), We consider
$\hat{K}_m^\varepsilon\subset[-2,2]^{\Gamma_m^\varepsilon}$ the
projection of $K_m^\varepsilon$ to the coordinates of
$\Gamma_m^\varepsilon$. We prove that $d_{\Gamma_m^\varepsilon}$
fragments $\hat{K}_m^\varepsilon$ (this is equivalent to say that
$d_{\Gamma_m^\varepsilon}$ fragments $K_m^\varepsilon$). We have a
continuous surjection $\phi:L\To \hat{K}_m^\varepsilon$ given by
$\phi(z) = (mz_{\delta(\gamma)})_{\gamma\in\Gamma_m^\varepsilon}$.
Observe that $d_{\Gamma_n^\varepsilon}(\phi(z),\phi(z')) =
m\sup_{\gamma\in\Gamma}|z_{\delta(\gamma)}-z'_{\delta(\gamma)}|
\leq m d_{\Delta}(z,z')$, so $$d_{\Gamma_n^\varepsilon}(x,x')\leq
m\inf\{d_\Delta(z,z') : \phi(z)=x, \phi(z')=x'\}.$$
 Since $d_\Delta$ fragments $L$, the
conclusion
follows from Lemma~\ref{quotientfrag}.\\

$(2)\Rightarrow (1)$: Let us call $\varepsilon_n=2^{-n}$. Without
loss of generality, we assume that
$\Gamma_m^\varepsilon\cap\Gamma_{m'}^\varepsilon =\emptyset$ for
$m\neq m'$. For $\gamma\in \Gamma$ and $n\in\mathbb{N}$, call
$m_{n}(\gamma)$ to be the only $m$ such that
$\gamma\in\Gamma_m^{\varepsilon_n}$.\\

We consider again $\hat{K}_m^\varepsilon\subset
[a,b]^{\Gamma_m^\varepsilon}$ the projection of $K_m^\varepsilon$
to the coordinates $\Gamma_m^\varepsilon$. The metric
$d_{\Gamma_m^\varepsilon}$ fragments
$\hat{K}_m^\varepsilon$ so this is an RN compactum.\\

Let $L\subset
\prod_{n\in\mathbb{N}}\prod_{m\in\mathbb{N}}K_m^{\varepsilon_n}$
be the set consisting of all $x$ for which there exists
$g(x)\in[a,b]^\Gamma$ such that

$$d_{\Gamma^{\varepsilon_n}_m}(g(x),p_{mn}(x))\leq\varepsilon_n$$
where $p_{mn}(x)\in K_m^{\varepsilon_n}$ is the coordinate of $x$ in the factor $K_m^{\varepsilon_n}$.\\

Notice that this element $g(x)$ is uniquely determined by $x$,
since $p_{m_n(\gamma)n}(x)_\gamma\rightarrow g(x)_\gamma$ as
$n\rightarrow\infty$.The space $L$ is compact and the map $g:L\To
[a,b]^\Gamma$ is continuous (these two facts are easily checked
considering net convergence). Moreover, $g(L)\supset K$ (this
follows from part (b) of condition (2) in the theorem). Since the
class RN is closed under the operations of taking countable
products and closed subspaces, $L$ is RN
compact and $K$ is a continuous image of an RN compactum.\qed\\

\section{A way to construct continuous images of RN
compacta}\label{construction}

We think now of a continuous image of an RN compactum as a compact
space $K$ that satisfies part (2) of Theorem~\ref{ciRN}. For
simplicity, we shall assume that $[a,b]=[0,1]$,
$\Gamma=\mathcal{N}$ and $u:\mathcal{N}\To\mathcal{N}$ is just the
identity map. We can think of Theorem~\ref{ciRN} as giving a
\emph{constructive} process to produce continuous images of RN
compacta in the
following way:\\

\begin{itemize}
\item [\textbf{Step 1.}] Begin with a family $\{K_s :
s\in\mathbb{N}^{<\omega}\}$ of compact subsets of
$[0,1]^\mathcal{N}$ such that
$d_{\mathcal{N}_s}$ fragments $K_s$.\\

\item [\textbf{Step 2.}] For every $s$ of length $n$, consider the
$\frac{1}{2^n}$-$\mathcal{N}_s$-enlargement of $K_s$, $$[K_s] =
\{x\in [0,1]^\mathcal{N} : \exists y\in K_s \
d_{\mathcal{N}_s}(x,y)\leq\frac{1}{2^n}\},$$ and some compact set
$L_s\subset [K_s]$.\\

\item [\textbf{Step 3.}] Finally, $K =
\bigcap_{s\in\mathbb{N}^{<\omega}}L_s$ is a
continuous image of an RN compactum.\\
\end{itemize}

This would be the general procedure but still it does not seem
that much \emph{constructive}: How can we get compacta $K_s$ for
step 1? And how can get the sets $L_s$ of step 2? Well, some
canonical choices can help us in this task:\\

For step 1 we can begin with an RN compactum $L\subset
[0,1]^\mathcal{N}$ which is fragmented by $d_{\mathcal{N}}$ and
make $K_s = L$ for all $s$. We know plenty of concrete examples of
such objects as we shall describe later. For step 2, simply
observe that the sets $[K_s]$ described above are themselves
compact, so we can take $L_s = [K_s]$, in our case

$$L_s = \{x\in [0,1]^\mathcal{N} : \exists y\in L :
d_{\mathcal{N}_s}(x,y)\leq\frac{1}{2^n}\}.$$

 After
this we shall call the resulting ciRN compactum $\tilde{L} =
\bigcap_{s\in\mathbb{N}^{<\omega}}L_s\subset [0,1]^\mathcal{N}$.
Notice that it only depends on the RN compactum $L\subset
[0,1]^\mathcal{N}$ that we took as starting point of our
construction. It follows from the proof of Theorem~\ref{ciRN} that
$\tilde{L}$ is indeed a continuous image of some closed subset of
the countable power $L^\mathbb{N}$.\\

A trivial example of compact space $L\subset [0,1]^\mathcal{N}$
which is fragmented by $d_\mathcal{N}$ is a scattered compactum.
Recall that a compactum is scattered if every nonempty subset
contains an isolated point. Scattered compacta are indeed
fragmented by any metric, even the discrete metric. Scattered
compacta, being totally disconnected, are typically found as
compact subsets of the Cantor cube $L\subset \{0,1\}^\mathcal{N}$
and in this case, it follows immediately that also
$\tilde{L}\subset \{0,1\}^\mathcal{N}$ is totally disconnected. It
is a known fact that continuous images of RN compacta (indeed all
qRN compacta) which are totally disconnected are RN compacta.
Thus, scattered compacta seem not to be good starting points for
our procedure if we are looking for candidates to be
counterexamples to the problem
of the continuous images.\\

As we mentioned in the previous paragraph, the problem of the
continuous images has a positive solution in the case when the
image is totally disconnected. So we focus rather on connected
compacta, where the problem becomes harder. We shall obtain our
connected RN compactum by taking the convex hull of a scattered
compactum:\\

\begin{prop}\label{convexhull}
Let $S\subset\{0,1\}^\mathcal{N}$ be a scattered compactum and let
$\overline{co}(S)\subset [0,1]^\mathcal{N}$ be the closure of its
convex hull
in $[0,1]^\mathcal{N}$, then $d_\mathcal{N}$ fragments $\overline{co}(S)$.\\
\end{prop}

We notice that it is not known in general whether the closed
convex hull of an RN compact space is again RN compact (when such an
operation makes sense and produces a compact set). This is indeed
a particular instance of the
problem of the continuous images, cf.~\cite{Namioka}.\\

Proof: Let $P(S)$ denote the space of regular Borel probability
measures on $S$ endowed with the weak$^\ast$ topology. Recall that
every continuous function $f:S\To [0,1]$ induces a continuous
function $\hat{f}:P(S)\To [0,1]$ given by $\hat{f}(\mu) = \int f
d\mu$. Applying this fact to every coordinate function over $S$,
we find a natural continuous function $g:P(S)\To
[0,1]^\mathcal{N}$ whose image is precisely $g(P(S)) =
\overline{co}(S)$. To get convinced about this latter fact, notice
that $g(\delta_s) = s$ for every $s\in S$ ($\delta_s$ denotes the
corresponding Dirac measure) and $g$ commutes with convex linear
combinations, and recall the well know fact that $P(S) =
\overline{co}\{\delta_s : s\in S\}$. Now, we can view $P(S)$ as
weak$^\ast$ compact subset of $C(S)^\ast$ where $C(S)$ is the
space of continuous functions over $S$. Since $S$ is scattered,
$C(S)$ is an Asplund space and this implies that any weak$^\ast$
compact subset of $C(S)^\ast$, and in particular $P(S)$, is fragmented by
the norm of $C(S)^\ast$ (cf. \cite{Fabianbook}). We observe the
following: for every $x,y\in P(S)$,
$$d_\mathcal{N}(g(\mu),g(\nu))\leq \|\mu-\nu\| =
\sup\{|h(\mu)-h(\nu)| : h\in C(S), \|h\|\leq 1\}.$$ Hence, the
pseudometric $d_\mathcal{N}(g(\mu),g(\nu))$ fragments $P(S)$ as
well, and making use of Lemma~\ref{quotientfrag}, we conclude
that $d_\mathcal{N}$ fragments $g(P(S))=\overline{co}(S)$.$\qed$\\

The map $g$ appearing in the proof of Proposition~\ref{convexhull}
is one-to-one provided that the linear span of the coordinate
functions is dense in $C(S)$, and in that case we would have that
$\overline{co}(S) = P(S)$. This will happen if for instance the
coordinate functions form an algebra of functions
over $S$, by the Stone-Weierstrass theorem.\\

Summarizing, the ciRN compacta that we are proposing are those
obtained in the following way:\\

\begin{itemize}
\item[\textbf{Step 0.}] Begin with a scattered compactum
$S\subset\{0,1\}^\mathcal{N}$.\\

\item[\textbf{Step 1.}] Consider its closed convex hull $L(S) =
\overline{co}(S)\subset [0,1]^\mathcal{N}$

\item[\textbf{Step 2.}] Finally, take the ciRN compactum
$$\tilde{L}(S) = \bigcap_{s\in\mathbb{N}^{<\omega}}L(S)_s =
\bigcap_{s\in\mathbb{N}^{<\omega}}\{x\in [0,1]^\mathcal{N} :
\exists y\in L(S)\  d_{\mathcal{N}_s}(x,y)\leq\frac{1}{2^n}\},$$
where $n$ is the length of $s$.
\end{itemize}

\begin{prob}
Let $S\subset \{0,1\}^\mathcal{N}$ be a compact scattered space.
Is the space $\tilde{L}(S)$ an RN compact?\\
\end{prob}

There are a couple of cases when we know that the answer to this
question is positive. One case occurs if the weight of $S$ is less
than $\mathfrak{b}$ \cite{cardinalb}. This will not interfere if
we take $S$ to be of weight the continuum, or at least weight
$\mathfrak{b}$, which on the other hand is the natural choice. The
other case is when $S$ is Eberlein compact, which implies that
$\tilde{L}(S)$ is Eberlein compact as well. We do not know much
more out of these
two cases, and those which can be obtained by mixing them up.\\

\section{The space of almost increasing functions}

We promised concrete examples and we are going to describe a very
concrete one in this section. The only variable on which the space
$\tilde{L}(S)$ depends is on the choice of some scattered compact
$S\subset \{0,1\}^\mathcal{N}$, which by the remarks above we
should take care not to be Eberlein compact. One natural
example of such is an ordinal interval $[0,\alpha]$ endowed with
the order topology. The way of viewing a scattered compactum of this
type as a subset of $\{0,1\}^\mathcal{N}$ is to fix a well order
$\prec$ on $\mathcal{N}$ and to declare

$$ S = \left\{x\in \{0,1\}^\mathcal{N} : \forall i\prec j\ x_i \leq x_j  \right\}$$

It is easy to check that the closed convex hull of this is nothing
else than

$$ \overline{co}(S) = \left\{x\in [0,1]^\mathcal{N} : \forall i\prec j\ x_i \leq x_j\right\} $$

It requires a little bit more work to realize that $\tilde{L}(S)$
has also a nice description as the set of ``almost increasing
functions''. We consider the following distance defined
on $\mathcal{N}$: $d(\sigma,\tau) = \frac{1}{2^{\min(n : \sigma_n\neq\tau_n)-1}}$.\\

\begin{thm}
In this case, we can describe the space $\tilde{L}(S)$ as follows:
$$\tilde{L}(S) = \left\{x\in [0,1]^\mathcal{N} : \forall \sigma\prec \tau\ x_\sigma \leq x_\tau + d(\sigma,\tau)\right\} $$
\end{thm}

Proof: It is enough to observe that for $x\in[0,1]^\mathcal{N}$
and $\varepsilon>0$, the following two conditions are
equivalent:\begin{enumerate} \item For every $\sigma\prec\tau$,
$x_\sigma\leq x_\tau+2\varepsilon$. \item There exists $y\in L(S)$
such that $|x_\sigma-y_\sigma|\leq\varepsilon$ for all $\sigma$.
\end{enumerate}
Namely, if condition $(1)$ holds, we can define $y_\sigma = \inf\{1,x_\tau+\varepsilon : \tau\succeq \sigma\}$.\qed\\

\section{RN quotients of a qRN compactum}

It is not only that we do not know whether every qRN compactum is
RN, moreover we do not know the answer to such a question as
whether every qRN compactum has an RN quotient of the same weight.
It is shown in~\cite{cardinalb} that every qRN compactum is a subspace of a
product of $\mathfrak d$ many RN compact spaces, and also that every qRN compactum of weight less than $\mathfrak b$ is RN compact ($\mathfrak d$ and $\mathfrak b$ denote the domination and bounding cardinal numbers, following the notation in \cite{vanDouwen}). As corollary one gets:\\

\begin{prop}
Let $K$ be a qRN compactum of weight $\kappa$.
\begin{enumerate}
\item If $\kappa>\mathfrak d$, then $K$ has an RN quotient of
weight $\lambda$ for every $\lambda<\kappa$. \item If
$cf(\kappa)>\mathfrak d$, then $K$
has an RN quotient of weight $\kappa$. \item Every quotient of $K$ of weight less than $\mathfrak{b}$ is RN compact.\\
\end{enumerate}
\end{prop}

Observe that the previous statements give no information about a compact space $K$ of weight $\mathfrak c$ under CH. We do not know whether the space of almost increasing functions is RN compact, but at least we can show the following.\\

\begin{thm}\label{almostincreasingquotient}
Let $K\subset [0,1]^\mathcal{N}$ be the compact space of almost
increasing functions associated to a well order of $\mathcal{N}$.
\begin{enumerate}\item For every cardinal $\lambda<\mathfrak c$, $K$
has an RN quotient of weight $\lambda$. \item If $\mathfrak c$ is
a regular cardinal, then $K$ has an RN quotient of weight
$\mathfrak c$. \item If the well order $(\prec)$ is chosen so that
every infinite $(\prec)$-interval is dense in the Baire space
$\mathcal{N}$,
then $K$ has an RN quotient of weight $\mathfrak c$.\\
\end{enumerate}
\end{thm}

Before entering the proof, we introduce some auxiliary definitions and results. For $\sigma\in \mathbb{N}^\mathbb{N}$, we call a
$\sigma$-good sequence a $\prec$-increasing sequence $\tau^\ast:
\tau^1\prec\tau^2\prec\cdots$ such that $d(\sigma,\tau^n)\leq
2^{-n-2}$. For $\tau^\ast$, $\sigma^\ast$ two sequences in
$\mathcal{N}$ we say $\tau^\ast\prec\sigma^\ast$ if
$\tau^i\prec\sigma^j$ for all $i,j$. A family $Z$ of such sequences is called separated if for any $\zeta^\ast,\xi^\ast\in Z$, either $\zeta^\ast\prec \xi^\ast$ or $\xi^\ast\prec \zeta^\ast$, so that $(\prec)$ is a well order on $Z$.\\

Given a $(\prec)$-increasing sequence $\tau^\ast$ in
$\mathcal{N}$ we produce a continuous function
$\phi[\tau^\ast]:K\To [0,1]$ as follows. First, we consider the
map $\Phi:[0,1]^\mathbb{N}\To [0,1]^\mathbb{N}$ which consists in
associating to each sequence $(x^n)$ a sequence $(y^n)$ with the
property that $|y^n-y^{n+1}|\leq 2^{-n}$, that we construct recursively: given $y^n$, we choose $y^{n+1}$ to be
closest number to $x^{n+1}$ that satisfies $|y^n-y^{n+1}|\leq 2^{-n}$; in a formula:
$$\Phi((x^k)_{k\in\mathbb{N}})_{n+1} = y^{n+1} = x^{n+1} + \text{sign}(y^n-x^{n+1})\cdot \max(|y^n-x^{n+1}|-2^{-n},0).$$

Notice that $\Phi$ is continuous because this recursive formula is continuous on $x_{n+1}$ and $y_n$, so $y^{n+1}$ depends continuously on $x_1,\ldots,x_{n+1}$.
For $x\in K\subset [0,1]^\mathcal{N}$, we also call
$\Phi[\tau^\ast](x) = \Phi((x_{\tau^n})_n)$. Clearly, the image of
$\Phi$ consists of convergent sequences so one can define
$$\phi[\tau^\ast](x) = \lim\Phi((x_{\tau^n})_{n\in\mathbb{N}})$$

The function $\phi[\tau^\ast]$ is continuous on $K$. The reason is that for every $m$, the map $x\mapsto \Phi[\tau^\ast](x)_n$ is continuous, and $\phi[\tau^\ast](x)$ is the uniform limit of this sequence of continuous functions, $|\Phi[\tau^\ast](x)_n - \Phi[\tau^\ast](x)_{n+1}|\leq 2^{-n}$.\\

\begin{lem}
Let $\sigma\in\mathcal{N}$ and let $Z$ be a separated family of
$\sigma$-good sequences. Then, $\mathcal{F}= \{\phi[\zeta^\ast] :
\zeta^\ast\in Z\}$
is a Namioka family.\\
\end{lem}

Proof: Given $m\in\mathbb{N}$ and $L$ a closed subset of $K$ we
will try to find a nonempty open set $V$ of $L$ of uniform
diameter less than
$2^{3-m}$ for the uniform metric associated to the family $\mathcal{F}$.\\

For every $n\in\mathbb{N}$, $\sigma|n =
(\sigma_0,\ldots,\sigma_{n-1})$, and we call $d_n =
d_{\mathcal{N}_{\sigma|n}}$ the corresponding pseudometric. Since
$d_n$ $\frac{1}{2^n}$-fragments $K$ we can find a nonempty open
set $U\subset L$ with $d_n$-diameter less than or equal to
$\frac{1}{2^n}$, for every $n\leq m$.
 On the other hand, $\zeta^n\in \mathcal{N}_{\sigma|n}$ for every $\zeta\in Z$ because $d(\zeta^n,\sigma)\leq 2^{-n-1}$. This means that for every
$\zeta\in Z$ and every $n\leq m$, the set $\{x_{\zeta^n} : x\in
U\}$ lies in an interval, say $I_{\zeta^n}$, of diameter less than
or equal to
$2^{-n}$.\\

\emph{Remark A}: We note that, if it would be the case that
$\Phi[\zeta^\ast](x)_n\in I_{\zeta^m}$ for every $x\in U$ and
every $\zeta\in Z$ then we would be done, because this would imply
that, for every $x\in U$, $\phi[\zeta^\ast](x)$ lies
in an interval of diameter at most $\frac{1}{2^{m-2}}$, namely $I_{\zeta^m}+\frac{1}{2^{m-1}}$.\\

If the the condition expressed in the remark above fails, we begin
a procedure of reducing $U$ and changing the intervals $I_n$ as
follows. If it is not the case it means that there exists $x\in U$
and some $\zeta$ such that $\Phi[\zeta^\ast](x)_m\not \in
I_{\zeta^m}$. We choose $\zeta^\ast_1$ to be the $(\prec)$-minimum
$\zeta^\ast$ to have this property and fix the corresponding $x\in
U$. In particular we will have that $\Phi[\zeta^\ast](x)_n\neq
x_{\zeta^n}$. By the definition of $\Phi$, if this happens it is
because there existed some $j<n$ such that $|x_{\zeta^j} -
x_{\zeta^{j+1}}|>2^{-j}$. Actually, since $x$ is almost increasing
and $d(\zeta^j,\zeta^{j+1})\leq 2^{-j-1}$ (for $\zeta^\ast$ is
$\sigma$-good), it must be the case that
$x_{\zeta^j}+2^{-j}<x_{\zeta^{j+1}}$.\\

 We define

\begin{eqnarray*} U_1 = \{ y\in U &:& |y_{\zeta^m}-x_{\zeta^m}|< 2^{-m-2}\\ &\text{ and }&
|y_{\zeta^j}-x_{\zeta^j}|< 2^{-m-5}\\ &\text{ and
}&|y_{\zeta^{j+1}}-x_{\zeta^{j+1}}|< 2^{-m-5}\} \end{eqnarray*}

\emph{Claim B}: For every $y\in U_1$ and every $\xi^\ast\succ \zeta^\ast$, we have that $x_{\zeta^j}+2^{-j-1} < y_{\xi^j}$.\\

Proof of the claim:
$$x_{\zeta^j}+2^{-j} < x_{\zeta^{j+1}} \leq y_{\zeta^{j+1}}+2^{-m-5}\leq y_{\xi^j} + d(\zeta^{j+1},\xi^j) + 2^{-m-5}\leq y_{\xi^j} + 2^{-j-1}.$$

Notice that for every $y\in U_1$:\\

\begin{itemize}
\item  If $\xi^\ast < \zeta^\ast$, $\Phi[\xi^\ast](y)_m\in I_m$,
\item  $|\Phi[\zeta^\ast](y)_m - \Phi[\zeta^\ast](x)_m|<2^{-m-2}$.\\
\end{itemize}

Therefore, if in addition $\Phi[\xi^\ast](y)_m\in I_{\xi^m}$ for every $\xi^\ast>\zeta^\ast$, the proof would be finished, on the same grounds as in the above \emph{Remark A}. If not, we must repeat the procedure of reduction and pass to a new further open set $U_2$, and look now at the minimum $\zeta^\ast_2>\zeta^\ast$ such that $\Phi[\zeta_2^\ast](x)_n\not \in I_{\zeta_2^n}$, etc. In order to conclude, we observe that the reduction procedure cannot be repeated infinitely many times, so that for some $k\in \mathbb{N}$, an open set $U_k$ of small $\mathcal{F}$ diameter will be found. The reason for the impossibility of infinite repetitions is precisely the inequality of \emph{Claim B}. There are only finitely many possible values for the natural number $j\leq m$, so for some $j$, after infinitely many steps we would have a sequence $x^{(1)}_{\zeta(1)} < x^{(2)}_{\zeta(2)} < \cdots$ of numbers in $[0,1]$ with $x^{(n)}_{\zeta(n)} + 2^{-j-1}< x^{(n+1)}_{\zeta(n+1)}$. This is a contradiction which finishes the proof of the lemma.$\qed$\\

\begin{lem}
Let $\sigma\in\mathcal{N}$, let $Z$ be a separated family of
$\sigma$-good sequences of cardinality $\kappa$, and $\mathcal{F}=
\{\phi[\zeta^\ast] : \zeta\in Z\}$. Consider $L$ the quotient
space of $K$ induced by the equivalence relation $(x\sim y \iff
\forall f\in\mathcal{F} f(x)=f(y))$. Then $L$ is an RN
compactum of weight $\kappa$.\\
\end{lem}

Proof: It is clear that $L$ is RN compact since $\mathcal{F}$
is a Namioka family. The point is to prove that the weight of $L$
is not less than $\kappa$. For this we shall show that actually
$L$ contains a copy of the compact space $P[0,\kappa]$ of
probability measures on the ordinal interval $[0,\kappa]$.
Remember that this space can be identified with

$$P[0,\kappa] = \{f:[0,\kappa]\To [0,1] : \alpha\prec\beta\Rightarrow f(\alpha)\leq f(\beta)\}.$$

Let $Z = \{\zeta^\ast(\alpha) : \alpha\prec \kappa\}$ be an
$(\prec)$-increasing enumeration of $Z$. To every $f\in
P[0,\kappa]$, we associate the function $\psi(f)\in K$ defined as\

\begin{itemize}
\item $\psi(f) (\sigma) = 0$ if $\sigma\prec \zeta^1(0)$ \item
$\psi(f)(\sigma) = f(\alpha)$ if $\zeta^1(\alpha)\preceq \sigma\prec
\zeta^1(\alpha+1)$
\item $\psi(f) (\sigma)=1$ if $\sigma\succ \zeta^1$ for all $\zeta\in Z$.\\
\end{itemize}

Then $\psi: P[0,\kappa]\To K$ is a one-to-one continuous map. The
key fact is that after composing with the quotient map onto $L$,
$$\psi:P[0,\kappa]\To K\To L$$
the function is still one-to-one. The reason is that
$\phi[\zeta^\ast(\alpha)](\psi(f)) = f(\alpha)$, for every
$\alpha\prec \kappa$ and every $f\in
P[0,\kappa]$.$\qed$\\

Proof of Theorem \ref{almostincreasingquotient}: In the case of item (3),
it is clear that, for any $\sigma\in\mathcal{N}$, we can find a
separated $\sigma$-good family of size $\mathfrak c$, since
actually we can find a $\sigma$-good sequence inside any infinite
$(\prec)$-interval. For (1) and (3), we shall show that for every
cardinal $\lambda\leq\mathfrak c$ there is a separated family of
$\sigma$-good sequences of cardinality $cf(\lambda)$. Let
$o:\mathcal{N}\To Ord$ be the function which associates to every
$\sigma\in\mathcal{N}$ the ordinal which indicates its position in
the well order $(\prec)$. Let us pick a point
$$\sigma\in\bigcap_{o(\tau)<\lambda}cl_\mathcal{N}\{ \xi :
o(\tau)\prec o(\xi)\prec \lambda\},$$ where
$cl_\mathcal{N}(\cdot)$ indicates the closure in the Baire space
of irrationals $\mathcal{N}$ with its usual metric topology. Such
a point $\sigma$ exists because $\mathcal{N}$ is herediatarily
Lindelof and we suppose that $cf(\lambda)>\omega$. It is clear
that we can find a separated family $Z$ of $\sigma$-good sequences
of size $cf(\lambda)$, all of them below the ordinal position
$\lambda$.\\

\end{document}